\documentclass[oneside]{amsart}
\usepackage{amsmath}
\usepackage{amsfonts}
\usepackage{amssymb}
\usepackage{latexsym}
\usepackage{amsthm}
\usepackage{setspace}
\usepackage{hyperref}
\usepackage{geometry}
\usepackage[usenames,dvipsnames]{color}

\everymath{\displaystyle}
\usepackage{url}
\newcommand{\Q}{\mathbb Q}
\newcommand{\A}{\mathbb A}

\DeclareMathOperator{\GL}{GL}

\DeclareMathOperator{\Gal}{Gal}

\DeclareMathOperator{\NCM}{\bf NCM}

\DeclareMathOperator{\LO}{{\bf LO}}

\newtheorem{thm}{Theorem}

\newtheorem{prop}[thm]{Proposition}

\newtheorem{qst}[thm]{Question}
\newtheorem{conj}[thm]{Conjecture}
\theoremstyle{remark}

\begin{document}

\title{Possible connection between a generalized Maeda's conjecture and local types}

\author{Luis Dieulefait}
\address{Facultat de Mathematiques, Universitat de Barcelona}
                     \email{ldieulefait@ub.edu}

\author{Panagiotis Tsaknias}
\address{Facult des Sciences, de la Technologie et de la Communication, Universit\'e du Luxembourg}
%                          Campus Kirchberg, 6 rue Richard Coudenhove-Kalergi, L-1359 Luxembourg}
\email{panagiotis.tsaknias@uni.lu}
\email{p.tsaknias@gmail.com}

%\date{May 2012}

\begin{abstract}
Here we follow on the proposed generalization of Maeda's conjecture made in \cite{Tsaknias2012a}. We report on computations that suggest a relation between the number of local types and the number of non-CM newform Galois orbits.
We extend the conjecture into spaces with non-trivial Nebentypus and provide a formula for the number of non-CM orbits for all levels and trivial Nebentypus. We also provide some numerical evidence towards further generalizations of this conjecture to totally real fields as well as further strengthening of it by proposing a structure for the corresponding Galois groups. 
\end{abstract}

\maketitle

%\tableofcontents

\section{Introduction}

A well-known conjecture of Maeda suggests that the number of newforms of level $1$ form a single orbit for all weights $k \geq 12$. Recently one of us proposed (in \cite{Tsaknias2012a}) a generalization of Maeda's  conjecture to arbitrary levels $N$. In particular, in \cite[Conjecture 2.2]{Tsaknias2012a}, the number of non-CM newform Galois orbits of fixed level $N$ and varying weight $k$ was conjectured to be eventually constant as a function of the weight $k$ (when we say "orbit" we are considering the action of the absolute Galois group of $\Q$ on the fields of coefficients of newforms). Furthermore the function $\NCM(N)$ giving this constant was conjectured to be multiplicative. This leads to the natural question: \emph{what can one say about $\NCM(p^n)$ for any prime $p$ and any integer $n\geq1$}?

As already mentioned in \cite{Tsaknias2012a}, Stein had proposed that the number of non-CM Galois orbits should be determined by the number of Atkin-Lehner decomposition factors. In the case of a prime $p$, this suggests the existence of two orbits: There is only one Atkin-Lehner operator and therefore only two possible factors one for each eigenvalue ($\pm1$). This was verified in \cite{Tsaknias2012a} for all primes $p\leq 200$: $\NCM(p)=2$ for all these primes. Moreover there was exactly one orbit for each eigenvalue.

Unfortunately, for $n\geq 2$, the Atkin-Lehner operators do not refine the new subspace enough in order to have each factor corresponding to exactly one orbit. Motivated by this problem, we propose the following:
\begin{conj}\label{conj:main}
The number of non-CM newform Galois orbits of level $p^n$, weight $k$ and trivial Nebentypus is eventually independent of the weight $k$ and equal to $\LO(p^n)$, the number of possible pairs $([\tau], \lambda_{AL})$, where $[\tau]$ is the Galois orbit of an inertial local type $\tau$ of conductor $p^n$ and $\lambda_{AL}$ an Atkin-Lehner eigenvalue compatible with $\tau$.%combinations with the same level and (class of) nebentypus.
\end{conj}

We have of course checked this in almost all computationally accessible cases, with no exceptions so far.

The fact that every possible inertial local type is actually occurring for all weights big enough (see \cite{Weinstein2009}) gives a lower bound on $\NCM(p^n)$,. The fact that it is multiplicative is in accordance with the conjecture formulated in \cite{Tsaknias2012a}. One thus has the following trivial observation:
\begin{prop}\label{prop:NCMlowbound}
Let $N=p_1^{e_1}\cdots p_s^{e_s}\geq1$, with $p_i$ distinct primes. Then
$$\NCM(N)\geq\LO(N),$$
where
$$\LO(N)=\prod_{i=1}^s\LO(p_i^{e_i})$$
\end{prop}
The Conjecture above essentially says that one has an equality in the previous proposition.  Given the existence result of Weinstein, this equality amount to unicity, i.e., the conjecture predicts that there is a unique Galois orbit of non-CM newforms of large weight $k$ with given local behavior at the primes in the level.

We also provide the following formula for $\LO(p^n)$ (we will include a proof of it in the next version of this preprint):
\begin{prop}\label{prop:LOformula}
Let $p$ be a prime number. If $p>2$ then
$$\LO(p^n)=
\begin{cases}
	1						&n=0							\\
	2						&n = 1							\\
	{\rm d}(p-1)+{\rm d}(p+1)-1	&n =2							\\
	{\rm d}(p-1)+{\rm d}(p+1)		&n=2m, m>1						\\
	4						&p>3\textrm{ and $n=2m+1$, or }n = 3	\\
	8						&p=3\textrm{ and }n =2m + 1> 3
\end{cases}$$
where ${\rm d}(a)$ is the number of divisors of $a$. For $p=2$ we have:
$$\LO(2^n)=
\begin{cases}
	1				&n=0 {\rm\ or\ } 2					\\
	2				&n = 1 {\rm\ or\ } 3					\\
%	1				&n =2							\\
%	2				&n=3							\\
	6				&n=4							\\
	4				&n=5							\\
	16				&n=6							\\
	8				&n\geq7,{\rm\ odd}				\\
	12				&n\geq8,{\rm\ even}				\\
%	8				&n=9							\\
%	12				&n=10							
\end{cases}$$
\end{prop}

We have gathered a substantial amount of computational data supporting the conjecture above. We have computed for every level up to $200$ and weight up to $40$ for almost all of them. The actual number of orbits converged to the conjectural one fairly quickly. We have also examined many levels up to $11^3$ that are prime powers but for a smaller number of weights, to check the validity of our formulas in this case. Even though our computational reach was smaller, we managed to observe convergence to the conjectured value fairly quickly (i.e. after weight $4$ or $6$).

We should stress at this point that there is nothing special with the $\Gamma_0$ level structure and that similar formulas can be derived in the $\Gamma_1$ case for any suitable choice of nebentypus, following exactly the same idea. We have also gathered data in support of this but not as extensive as in the trivial nebentypus case.

Let now $\Q_f$ be the coefficient field of a newform $f$ and $\Q^{\textrm{gal}}_f$ a Galois closure of it. The original conjecture from Maeda was actually a statement on the Galois group $Gal(\Q^{\textrm{gal}}_f/\Q)$: The claim is that in the case of level $N=1$ this group is isomorphic to $S_[\Q_f:Q]$, the symmetric group on $[\Q_f:\Q]$ elements. Motivated by this and following the same principle, i.e. that the only reason the conjecture for general $N$ is deviating from Maeda's original formulation should be the obvious ones. Let's assume for simplicity that $f$ has no inner twists. Then it is well-known, and can be easily seen\footnote{For example the fixed field of permuting just the inertial types in an orbit}, that $\Q_f$ has an abelian subfield $L_f$ (this field in some cases is just $\Q$, for example in the square-free level case). This forces the group $\Gal(\Q^{\textrm{gal}}_f/\Q)$ to have $\Gal(L_f/\Q)$ as a quotient. We pose the following question:
\begin{qst}
Let $f$ be a non CM newform, $\Q_f$ its coefficient field and denote by $\Q^{\textrm{gal}}_f$ and $L_f$ the Galois closure and the obvious abelian subfield of $\Q_f$ respectively. Is  the following exact sequence true: 
$$1 \to S_[\Q_f:L_f] \to \Gal(\Q^{\textrm{gal}}_f/\Q) \to \Gal(L_f/\Q) \to 1$$
\end{qst}

The amount of data gathered towards this question (\cite{Calic2015}) are perhaps not sufficient but we are tempted to answer  yes for high enough $k$. In the case of inner twists, there are cases where an extra term has to be added at the beginning of the above sequence. We hope to get into more details in the upcoming and more detailed report. Notice however that their presence is automatically excluded in the case of square-free level and trivial nebentypus.

This we believe provides an adequate generalization of Maeda's conjecture to arbitrary level and nebentypus in the case of classical modular forms. We felt tempted to try this approach on Hilbert Modular forms over arbitrary totally real fields $F$. Following the idea that obvious exceptions had to be considered we immediately excluded the CM forms from a possible statement like in the classical case, but it became immediately obvious (e.g by looking at trivial level) that one has to exclude Base Change forms as well. The computational data we have gathered point to the same direction as the Conjecture above: After excluding the aforementioned exceptional cases, the number of newform Galois orbits over a totally real field $F$ for a given level $\mathfrak{n}$ is eventually constant as a function of the weight $k$. Moreover this constant is equal\footnote{It has been suggested to us by Dembele that for this one needs to consider fields $F$ with trivial narrow class number.} to the number of all possible combinations of pairs of compatible inertial types and Atkin-Lehner eigenvalues. Again, since an existence result (of a Hilbert modular forms with given local behavior at the primes in the level) is known for sufficiently large weight, this conjecture amounts to unicity of the Galois orbit of Hilbert newforms with given local behavior (for weight sufficiently large). In a forthcoming work of the first author and Ariel Pacetti, a proof of arbitrary base change (i.e., to an arbitrary totally real extension of the base field) for Hilbert modular forms will be given under the assumption of this generalization of Maeda's conjecture to the Hilbert case.

At this point we believe it is fairly natural to try to further extend this sort reasoning: If we view Hilbert modular forms over $F$ as automorphic forms of $\GL_2(\A_F)$, we feel it is tempting to formulate and test similar statements over other algebraic groups such as $\GL_n$, or the ones associated with units of quaternion algebras. In the later case for example there are cases were one immediately gets a Maeda like conjecture by utilizing the Jacquet-Langlands correspondence: Consider the algebraic group $G$ associated with the quaternion algebra over $\Q$ ramified at two distinct primes $p,q$. Then the automorphic forms of $G$ of trivial level structure are in correspondence with the classical newforms over $\Q$ of level $p, q$ and therefore we expect them to group into $4$ orbits as Conjecture \ref{conj:main} suggests. This list is by no means exhaustive and we must admit that we have gathered no computational evidence towards these kind of generalizations but it is our hope it will spark some research interest towards this direction.

Finally, it is perhaps interesting to contrast this conjecture with a Question posed by Buzzard and which proposes a completely different behavior for the $G_p$-Galois orbits, where $p$ is any rational prime: The degree of the coefficient fields as an extension of $\Q_p$ for fixed level $N$, fixed prime $p$ and varying weight $k$ is bounded by a constant depending on
$N$ and $p$.


\begin{thebibliography}{9}

\bibitem{Calic2015}
L.~Calic and V.~Scansi.
\newblock Galois groups in generalisations of Maeda's conjecture.
\newblock {\em Master Thesis}, 2015

\bibitem{Tsaknias2012a}
P.~Tsaknias.
\newblock A possible generalization of Maeda's conjecture.
\newblock {\em Computations with modular forms}, Springer, Cham, 317--329, 2014

\bibitem{Weinstein2009}
J.~Weinstein.
\newblock Hilbert Modular Forms with Prescribed Ramification.
\newblock {\em International Mathematics Research Notices}, 2009(8):1388--1420, 2009.

\end{thebibliography}
\end{document}